\documentclass[11pt,reqno]{article}
\usepackage{graphicx,amsmath,amssymb,setspace}
\usepackage{natbib}

\newtheorem{thm}{Theorem}[section]
\newtheorem{lemma}{Lemma}[section]
\newtheorem{remark}{Remark}[section]

\newenvironment{proof}{\textsc{Proof:}}{\mbox{ } \hfill $\Box$ \vspace{2mm}}

\numberwithin{equation}{section}

\newcommand{\E}{\mathbb{E}}
\newcommand{\R}{\mathbb{R}}
\renewcommand{\P}{\mathbb{P}}

\newcommand{\ls}{\lambda^S}
\newcommand{\lo}{\lambda^{O}}
\newcommand{\noi}{\noindent}
\renewcommand{\a}{\bar{a}}
\renewcommand{\L}{\mathcal{L}}
\newcommand{\A}{\mathbb{A}}
\newcommand{\Ac}{\mathcal{A}}

\title{\LARGE \bf
\sc Optimal Deferred Life Annuities to Minimize the Probability of
Lifetime Ruin
\thanks{We thank S. David Promislow for helpful
conversations concerning this project.} 
}

\author{Erhan Bayraktar 
\thanks{Department of Mathematics, University of Michigan, Ann Arbor, MI 48109,  email: erhan@umich.edu}
\thanks{E. Bayraktar is supported in part by the National Science Foundation under grant DMS-0604491.}
\and Virginia R. Young \thanks{ Department of Mathematics,
University of Michigan, Ann Arbor, Michigan, 48109,
email:vryoung@umich.edu} }
\date{}

\begin{document}
\maketitle

\begin{abstract}
We find the minimum probability of lifetime ruin of an investor who can invest in a market with a risky and a riskless asset and who can purchase a deferred annuity. Although we let the admissible set of strategies of annuity purchasing process to be the set of increasing adapted processes, we find that the individual will not buy a deferred life annuity unless she can cover all her consumption via the annuity and have enough wealth left over to sustain her until the end of the deferral period.

\end{abstract}

\section{Introduction}

In this paper, we study the problem of investing in a risky financial market and buying deferred life annuities to minimize the probability of wealth reaching zero before dying, the so-called {\it probability of lifetime ruin}.  Although life annuities provide income security in retirement, very few retirees choose a life annuity over a lump sum.  For example, in a recent comprehensive Health and Retirement Survey, only 1.57\% of the respondents reported annuity income; similarly, only 8.0\% of respondents with a defined contribution pension plan selected an annuity payout. For a further discussion on the thinness of the annuity market see \cite{my} and the references therein.

In a well-cited paper from the public economics literature, \cite{yaari} proved that in the absence of bequest motives -- and in a deterministic financial economy -- consumers will annuitize all of their liquid wealth. \cite{richard} generalized this result to a stochastic environment, and recently \cite{davidoff} demonstrated the robustness of Yaari's result.  The common theme of this research is the assumption of a rational expected utility-maximizing economic agent.  While this von-Neumann-Morgenstern framework is the basis of most investment and consumption strategies studied in mathematical finance, it is perhaps difficult to apply as a tool for normative advice because of the subjectiveness of the individual's utility function.

Recently, a variety of papers in the risk and portfolio management literature have reintroduced the \cite{roy} Safety-First rule and applied the concept to probability maximization of achieving certain investment goals. For example, \cite{browne95} derived the optimal dynamic strategy for a portfolio manager who is interested in minimizing the probability of shortfall.  \cite{milevsky} introduced the probability of lifetime ruin as a risk-metric for retirees, albeit in a static environment.  As an extension of that work, \cite{young04} determined the optimal dynamic investment policy for an individual who consumes at a specific rate, who invests in a complete financial market, and who does not buy annuities.   In work related to this paper, \cite{mmy} extended \cite{young04} by allowing the individual to buy immediate life annuities.  They showed that the individual will not buy {\it any} amount of annuity income  unless she can cover {\it all} of her consumption.  In this paper, we extend \cite{young04} by allowing the individual to buy deferred life annuities.

The layout of the paper is as follows.  In Section 2, we present the financial market in which the individual can invest her wealth.  We also present a motivating example to demonstrate the individual's myopia when presented with the opportunity to buy an immediate life annuity.  In addition to investing in a risky asset, in Section 3, we allow the individual to purchase deferred life annuities but not immediate annuities.  The individual's goal is to minimize her probability of lifetime ruin, and we prove a verification lemma for this minimal probability.  It turns out that the individual exhibits the same myopia as witnessed in \cite{mmy} and in the simple example in Section 2.  Namely, she will not purchase any deferred annuities until she can purchase an annuity that will cover all her consumption with enough wealth remaining so that she will not ruin before the end of the deferral period.

In Section 4, we present a related optimal stopping problem.  We show that the Legendre transform of the value function of that stopping problem equals the minimal probability of lifetime ruin.  Because the value function of the optimal stopping problem solves a {\it linear} variational inequality, one can use methods from pricing American options to determine that value function.  Then, one can take its Legendre transform to obtain the minimal probability of lifetime ruin.  We demonstrate this method via a numerical example in Section 5.  We conclude in Section 6.

\section{Minimizing the Probability of Lifetime Ruin}

In this section, we describe the financial market in which the individual can invest her wealth, and we formulate the problem of minimizing the probability of lifetime ruin in this market.  Then, we allow the individual to purchase an immediate life annuity at {\it one} point in time.  We show that the individual will {\it not} do so unless she has enough money to buy an annuity that will cover {\it all} her consumption, even if she is not allowed to buy annuities at any time thereafter.

\subsection{Financial model}

We consider an individual aged with future lifetime described by the random variable $\tau_d$.  Suppose $\tau_d$ is an exponential random variable with parameter $\ls$, also referred to as the force of mortality or hazard rate; in particular, $\E[\tau_d] = 1/\ls$.  The superscript $S$ indicates that the parameter equals the individual's subjective belief as to the value of her hazard rate.

We assume that the individual consumes wealth at a constant {\it net} rate of $c$; this rate might be given in real or nominal units.  We say that the rate $c$ is a net rate because it is her rate of consumption offset by any current income she receives. One can interpret $c$ as the minimum net consumption level below
which the individual cannot (or will not) reduce her consumption further; therefore, the minimum probability of lifetime ruin that we compute gives a lower bound for the probability of ruin under any consumption function bounded below by $c$.

The individual can invest in a riskless asset, which earns interest at the rate $r \ge 0$.  Also, she can invest in a risky asset whose price satisfies
\begin{equation}
dS_t = \mu S_t dt + \sigma S_t dB_t, \quad S_0 = S > 0,
\end{equation}
\noi in which $\mu > r$, $\sigma > 0$, and $B$ is a standard Brownian motion with respect to a filtration $\mathbb{F} = \{{\cal F}_t \}$ of a probability space $(\Omega, {\cal F}, \P)$.  We assume that $B$ is independent of $\tau_d$, the random time of death of the individual.  If $c$ is given as a real rate of consumption (that is, after inflation), then we express $r$ and $\mu$ as real rates.

Let $\pi_t$ denote the amount invested in the risky asset at time $t$, and let $\pi$ denote the investment strategy  $\{\pi_t\}_{t \geq 0}$.  We say that a strategy $\pi$ is {\it admissible} if the process $\pi$ is adapted to the filtration $\mathbb{F}$ and if $\pi$ satisfies the condition $\int_0^t \pi_s^2 \, ds < \infty$, almost surely, for all $t \ge 0$.  The wealth dynamics of the individual for a given admissible strategy $\pi$ are given by

\begin{equation}\label{eq:wealth}
dW_t = [r W_t + (\mu - r) \pi_t - c] dt + \sigma \pi_t dB_t, \quad
W_0 = w \ge 0.
\end{equation}

By ``lifetime ruin," we mean that the individual's wealth reaches the value 0 before she dies.  Instead, one could define ruin as reaching some non-zero value $b$ before dying, in which $b$ could represent assets below which an individual qualifies for social assistance or considers herself impoverished.  Alternatively, $b$ could represent a minimum bequest amount that the individual wants to leave her heirs.  For simplicity, we take the ruin level to be $b = 0$.  Note that $0$ is, then, the absorbing boundary of the wealth process $W$ so that if wealth reaches $0,$ then the wealth process is ``killed'' and the game is over.

In this simple (time-homogeneous) setting, denote the minimum probability that the individual outlives her wealth by $\phi(w; c)$, given that the individual is alive, in which one minimizes the probability of ruin over admissible investment strategies. The argument $w$ indicates that one conditions the ruin probability on the individual possessing wealth $w$ the current time, and we explicitly represent the rate of consumption $c$ in this ruin probability.  \cite{young04} explicitly determines that $\phi$ is given by

\begin{equation} \label{eq:probruin_no}
\phi(w; c)=\left ( 1 - {rw \over c} \right )^d \; \text{ for } 0 \le w \le {c \over r},  \text{ in which}
\end{equation}
\begin{equation}\label{2.4}
d = {1 \over 2r} \left[ (r+\ls + m) + \sqrt{(r+\ls +m)^2-4r \ls} \; \right]>1 \; \text{and}
\end{equation}
\begin{equation}
m={1 \over 2}\left ( {\mu-r \over \sigma} \right )^2.
\end{equation}

Also, \cite{young04} shows that the optimal investment strategy $\{ \pi^* \}$ is given in feedback form.  Specifically, $\pi^*_t = \pi^*(W^*_t)$, in which $W^*$ is the optimally-controlled wealth, and in which $\pi^*$ is given by

\begin{equation}
\pi^*(w; c) ={\mu-r \over \sigma^2}\,{c - rw \over (d-1)r}  \; \text{ for } 0 \le w \le {c \over r}.
\end{equation}

\noi For $w \ge c/r$, any investment strategy is optimal if at least $c/r$ is invested in the riskless asset, so that at least $c$ is produced as continuous income risklessly.

\subsection{One-time purchase of an immediate life annuity}

\cite{mmy} show that if one introduces immediate life annuities into the model given in Section 2.1, then the individual will not buy an immediate life annuity until her wealth is large enough to cover all her consumption.  We demonstrate a similar phenomenon in this section.

Suppose that an individual can buy an immediate life annuity at {\it one} point in time, say now.  Assume that the price of a life annuity that pays \$1 per year continuously is given by

\begin{equation}\label{imm-ann}
\int_0^{\infty} e^{-r s} e^{-\lo s} \, ds = \frac{1}{\rho},
\end{equation}

\noi in which $\rho \triangleq r+\lo$, and $\lo > 0$ is the constant objective hazard rate that is used to price annuities.

Suppose $w \ge c/\rho$; then, it is optimal for the individual to spend $c/\rho$ to buy an immediate annuity that will pay at the continuous rate $c$ for the rest of her life.  In this case, the individual will not ruin, under the convention that if her net consumption rate becomes 0, then she is not considered ruined even if her wealth is 0.  (The latter occurs if $w = c/\rho$ immediately before buying the annuity.)

Next, suppose $w < c/\rho$, and suppose the individual spends $\Delta c/\rho$ on an immediate annuity so that her net consumption become $c - \Delta c$.  After buying the annuity, the probability of ruin equals $\phi(w - \Delta c/\rho; c - \Delta c)$, in which $\phi$ is given by (\ref{eq:probruin_no}).  Because $\phi$ is given explicitly, one can use elementary calculus to show that  $\phi(w - \Delta c/\rho; c - \Delta c)$ is strictly increasing with respect to $\Delta c$ for $\Delta c \in (0, w \rho)$.  Therefore, the optimal amount of immediate annuity income for the individual to buy is 0.  In other words, if the individual cannot completely eliminate her probability of ruin, then she will not buy any amount of immediate annuity.  She does not want to relinquish the ability to trade flexibly between the riskless and risky assets, even though the riskless return $r < \rho = r + \lo$.

Note that the individual exhibits surprising myopia in her decision not to buy an immediate annuity.  It is surprising because the life annuity offers the riskless rate of $\rho > r$. After this one-time opportunity to buy an immediate life annuity, she will only be able to lock into the riskless rate $r < \rho$.  If her wealth is less than $c/\rho$, then she gives up this chance to buy an immediate annuity  because by doing so, she would increase her probability of lifetime ruin.

In what follows, we observe a similar behavior when it comes to buying deferred life annuities.  Namely, the individual will not buy a deferred life annuity unless she can cover all her consumption via the annuity and have enough wealth left over to sustain her risklessly until the end of the deferral period.

\section{Deferred Life Annuities}

In this section and in the remainder of the paper, we consider the problem of minimizing the probability of lifetime ruin when the individual can purchase deferred annuities, in addition to investing in riskless and risky assets.

\subsection{Probability of lifetime ruin}

At time $t$, the individual can buy a deferred life annuity that pays at a continuous rate and begins at a fixed time $T > t$, if she is alive then, and that continues until she dies. Here one might take $T$ as the retirement time of the individual. The value of this deferred life annuity that pays \$1 per year continuously is given by

\begin{equation}
\a(t) = \int_T^{\infty}e^{-r(s-t)}e^{-\lo(s - t)}ds = \frac{1}{\rho}e^{-\rho (T - t)},
\end{equation}

\noi in which $\rho = r+\lo$, as in (\ref{imm-ann}).  Before and after time $T$, the individual can invest in the riskless and risky assets described in Section 2.

Let $A_s$ denote the cumulative amount of deferred annuity income that will begin at time $T$, purchased at or before time $s$, and let $\pi_s$ denote the amount invested in risky asset at time $s$, as in Section 2.  Let $\Ac$ denote the annuity-purchasing strategy $\{A_s\}_{s \in [t,T)}$.  We say that a strategy $(\Ac,\pi)$ is {\it admissible} if the non-decreasing process $\Ac$ is adapted to the filtration $\mathbb{F}$; if $A_s \in [0, c]$ for $s \in [t,T)$; and if $\pi$ is admissible as defined in Section 2.

Note that it is not optimal for an individual to spend wealth now to cover {\it more} than the consumption rate $c$.  Also, if the individual were to be in the fortunate position of receiving $A \ge c$ of deferred income at time $T$ either through prior purchases of deferred annuities or through pension income, then the effect would be as if $A$ were equal to $c$.  It follows that it is not a restriction to require that $\Ac$ satisfy $0 \le A_t \le c$ for all $t \in [0, T)$.

Denote the set of admissible strategies by $\A$.  The wealth dynamics of the individual for a given strategy $(\Ac,\pi)\in \A$ are given by

\begin{equation}\label{eq:wealthDA}
dW_s=
\begin{cases}
[rW_{s-} + (\mu-r) \pi_{s-} - c]ds + \sigma \pi_{s-} dB_s - \a(s) dA_s, & t<s<T,
\\ [rW_s + (\mu-r)\pi_s + A_T - c]ds+\sigma \pi_s dB_s, & s > T,
\end{cases}
\end{equation}

\noi where  $W_t = w \ge 0, \, A_t = A \in [0,c].$

In this setting, denote the minimum probability that the individual outlives her wealth by $\psi(w, A, t)$, given that the individual is alive at time $t$.  The arguments $w$, $A$, and $t$ indicate that one conditions the ruin probability on the individual possessing wealth $w$ at time $t$ and on the individual receiving deferred annuity income $A$ at time $T$ due to purchases at or before time $t < T$.  Note that we also allow for the individual
to possess deferred income in the form of pension or Social Security benefits; such amounts are implicit in $A$.

Let $\tau_0$ denote the first time that wealth equals 0, and recall that $\tau_d$ denotes the random time of death of the individual. Thus, $\psi$ is the minimum probability that $\tau_0 < \tau_d$, in which one minimizes with respect to admissible investment strategies $\pi$ dynamically and with respect to the deferred annuity benefits $\Ac$. The minimum probability of ruin can be written as

\begin{equation} \label{eq:probruin}
\psi(w, A, t) = \inf_{(\Ac,\pi)\in \A} \P \left[ \tau_0 < \tau_d \,\big| A_t=A, W_t=w \right],
\end{equation}

\noi $w \geq 0$, $A \in [0, c]$, and $t \geq 0$.

\begin{remark}
Suppose the individual purchased deferred annuities at or before time $t$ to the extent that her income from those annuities will be $A$ beginning at time $T$. If she buys another deferred annuity at time $t$ that covers all of the excess consumption $c-A$ with enough left over to live on until time $T$ by investing only in the riskless asset, she can avoid ruin for sure if she is wealthy enough at time $t$.  Next, we determine this so-called ``safe
level.''


If the individual begins with $W_{t-} = w$ wealth, buys $c - A$ of deferred annuity income, consumes at a rate of $c$, invests all her wealth in the riskless asset at rate $r$, then for $t < s < T,$ her wealth follows

\begin{equation}
dW_s=(rW_s - c)ds, \quad W_{t+} = w - {c - A \over \rho} e^{-\rho(T - t)}.
\end{equation}

\noi  It follows that

\begin{equation}
W_T=\left(w-\frac{c-A}{\rho}e^{-\rho (T - t)}\right)e^{r(T - t)}-c\frac{e^{r(T - t)} - 1}{r}.
\end{equation}

The safe level $w=\bar{w}(A,t)$ can be found by solving $W_T= 0,$ which yields

\begin{equation}\label{eq:bar-w}
\bar w(A,t)=c\frac{1-e^{-r(T - t)}}{r}+(c-A)\frac{e^{-\rho (T - t)}}{\rho}, \quad t \leq T, \, A \in[0,c].
\end{equation}

\noi Note that we assume that if total rate of income at time $T$ equals the rate of consumption, then the ruin probability is $0,$ even if wealth is $0$ at that time.

\end{remark}

In the following remark, we provide a useful characterization of the minimum probability of lifetime ruin in our setting.

\begin{remark}\label{rem:rep-of-prob-ruin}
Recall that $\tau_0 = \inf \{t > 0: W_t = 0 \}$. Then, because $\tau_0$ and $\tau_d$ are independent, we have the following expression for $\psi :$

\begin{equation}
\begin{split}
\psi(w, A, t) &= \inf_{(\Ac, \pi) \in \A} \E \left[ \int_t^\infty \ls e^{\ls (s - t)} \, {\bf 1}_{\{ t \le \tau_0 \le s \}} \, ds \Big| W_t = w, A_t = A \right]  \\
&= \inf_{(\Ac, \pi) \in \A} \E \left[ \int_{\tau_0}^\infty \ls e^{\ls (s - t)} \, {\bf 1}_{\{ t \le  \tau_0 < \infty \}} \, ds \Big| W_t = w, A_t = A \right]  \\
&= \inf_{(\Ac, \pi) \in \A} \E \left[e^{-\ls (\tau_0 - t)} \,{\bf 1}_{\{t \le \tau_0 < \infty\}} \Big| W_t = w, A_t = A \right] \\
& = \inf_{(\Ac, \pi) \in \A} \E \left[ e^{-\ls (\tau_0 - t)} \, {\bf 1}_{\{t \le \tau_0 \le T\}} + e^{-\ls (T - t)} \phi(W_T; c- A_T) \, {\bf 1}_{\{\tau_0 > T \}} \Big| W_t = w, A_t = A \right],
\end{split}
\end{equation}

\noi in which $\phi$ is the minimum probability of ruin given in
$(\ref{eq:probruin_no})$.

\end{remark}

\subsection{Verification lemma}

In this section, we provide a useful inequality which allows us in Section 4 to show that if we find a smooth solution to a given boundary-value problem, then that solution is the minimum probability of ruin defined in (\ref{eq:probruin}).

\begin{lemma}\label{lem:verf-lemma}
For any $\pi \in \mathbb{R},$ define the functional operators $\mathcal{L}^{\pi}$ and $\mathcal{L}$ through their actions on a test function $f$ as

\begin{equation}
\L^{\pi} f= f_t+ (rw-c)f_w+ \frac{1}{2}\sigma^2 \pi^2 f_{ww}-\ls f,
\end{equation}

\noi and

\begin{equation}
\L f \triangleq \min_{\pi} \L^{\pi}f= f_t+(rw -c)f_w- m
\frac{f_w^2}{f_{ww}}-\ls f.
\end{equation}

\noi Note that the second equality is satisfied if $f_{ww} > 0$. Let $v = v(w,A,t)$ be a non-increasing convex function of $w$ that is twice-differentiable with respect to $w,$ except possibly at $w = \bar w(A, t)$ where we assume that it has right- and left-derivatives, and differentiable with respect to $A$ and $t$. Suppose $v$ satisfies the following conditions:

\begin{enumerate}
\item $\L^{\pi}v (w,A,s) \geq 0,$
\item $\a(s)v_w(w,A,s)-v_A (w,A,s) \leq 0,$
\item $v(0, A, s) = 1,$
\end{enumerate}
for any $\pi \in \mathbb{R}$, $A \in [0,c],$ $w \geq 0,$ and $t \le s < T,$ and
\begin{enumerate}
\item[4.] $v(w,A,s)=\psi(w,A,s) = \phi(w; c - A),$ $s \geq T$.
\end{enumerate}
\noi Then,
\begin{equation}\label{eq:first-claim}
v(w,A,t) \leq \psi(w,A,t),
\end{equation}
\noi for all $w \geq 0$ and $A \in [0,c]$.
\end{lemma}

\begin{proof}
Let $\tau_n \triangleq \{s\geq t: \int_t^{s} \pi_s^2 ds \geq n\}$. Define $\tau \triangleq \tau_0 \wedge \tau_n$, which is a stopping time with respect to the filtration $\mathbb{F}$; then, by using It\^{o}'s formula for semi-martingales, we can write
\begin{equation}\label{eq:Ito}
\begin{split}
&e^{-\ls ((\tau \wedge T) - t)} v(W_{\tau \wedge T}, A_{\tau \wedge T},\tau \wedge T)=v(W_t,A_t,t)+\int_t^{\tau \wedge T}e^{-\ls (s-t)}v_w(W_s,A_s,s)\sigma \pi_s dB_s \\
&+\int_t^{\tau \wedge T}e^{-\ls (s-t)} \L^{\pi_s}v(W_s,A_s,s)ds-\int_t^{\tau \wedge T}e^{-\ls (s-t)}v_{w}(W_s,A_s,s)\a(s)d A_s^{(c)} \\
& +\int_t^{\tau \wedge T}e^{-\ls (s-t)}v_{A}(W_s,A_s,s)dA_s^{(c)}+\sum_{t \leq s \leq \tau \wedge T} e^{-\ls (s-t)} \left(v(W_s,A_s,s)-v(W_{s-},A_{s-},s-)\right).
\end{split}
\end{equation}
\noi Here, $\Ac^{(c)}$ is the continuous part of $\Ac$, that is,
\begin{equation}
A^{(c)}_t \triangleq A_t-\sum_{0 \leq s \leq t}(A_s-A_{s-}).
\end{equation}

Since $v$ is convex in $w$, $v^2_w(w,A,t) \leq v^2_w(0, A, t)$ for $w \ge 0$. Also, any pair $(A,t)$ belongs to the compact set $ [0,c] \times [0,T]$. Therefore,
\begin{equation}
\E\left\{\int_t^{\tau \wedge T}e^{-2 \ls (s-t)} v^2_{w}(W_s,A_s,s) \sigma^2 \pi_s^2 ds\bigg|W_t=w,A_t=A\right\}<\infty,
\end{equation}
\noi which implies that
\begin{equation}\label{eq:sint}
\E\left\{\int_t^{\tau \wedge T} e^{-\ls (s-t)} v_w(W_s,A_s,s)\sigma \pi_s dB_s\bigg|W_t=w,A_t=A\right\}=0.
\end{equation}
By taking expectations of equation (\ref{eq:Ito}), while using (\ref{eq:sint}) and Assumptions 1 and 2 in the statement of the verification lemma, we obtain
\begin{equation}\label{inequ}
\E\left\{e^{-\ls ((\tau \wedge T) - t)}v(W_{\tau \wedge T},A_{\tau \wedge T}, \tau \wedge T)\Big|W_t=w, A_t=A\right\} \geq v(w,A,t).
\end{equation}
\noi In deriving (\ref{inequ}), we also used the fact that
\begin{equation}
\sum_{t \leq s \leq \tau \wedge T} e^{-\ls (s-t)} \left(v(W_s,A_s,s)-v(W_{s-},A_{s-},s-)\right) \geq 0,
\end{equation}

\noi because Assumption 2 implies that $v$ is non-decreasing in the direction of jumps of the state process $\{W_s,A_s\}_{s \geq 0}$.

An application of the dominated convergence theorem to (\ref{inequ}) yields

\begin{equation}\label{eq:v-T-v-t}
\E\left\{e^{-\ls ((\tau_0 \wedge T) - t)}v(W_{\tau_0 \wedge T},A_{\tau_0 \wedge T}, \tau_0 \wedge T)\Big|W_t = w, A_t=A\right\} \geq
v(w,A,t).
\end{equation}

\noi By using Assumptions 3 and 4, one can rewrite (\ref{eq:v-T-v-t}) as

\begin{equation}\label{eq:vleqanyprob}
\begin{split}
v(w,A,t) &\leq \E\left\{e^{-\ls (\tau_0 - t)}v(\tau_0,A_{\tau_0},\tau_0){\bf 1}_{\{t \le \tau_0 \leq T \}}+e^{-\ls (T-t)} v(W_T,A_T,T){\bf 1}_{\{\tau_0 > T\}}\Big|W_t=w, A_t =A\right\} \\
&=\E\left\{e^{-\ls (\tau_0 - t)}{\bf 1}_{\{t  \le \tau_0 \leq T \}}+e^{-\ls (T - t)}\phi(W_T; c - A_T){\bf 1}_{\{\tau_0 > T\}}\Big|W_t=w, A_t=A\right\}.
\end{split}
\end{equation}
\noi This equation implies that
\begin{equation}\label{eq:vleq}
\begin{split}
v(w,A,t) &\leq \inf_{(\Ac,\pi)\in \A}\E\left\{e^{-\ls (\tau_0 - t)}{\bf 1}_{\{t \le \tau_0 \leq T \}}+e^{-\ls (T-t)} \phi(W_T; c -A_T){\bf 1}_{\{\tau_0 > T\}}\bigg|W_t=w, A_t=A\right\} \\
& =\psi(w,A,t),
\end{split}
\end{equation}
\noi in which the last equality follows from Remark~\ref{rem:rep-of-prob-ruin}. Thus, we have proven (\ref{eq:first-claim}).
\end{proof}

In the next section, we show that the minimum probability of ruin $\psi$ is intimately related to the solution of an optimal stopping problem.

\section{Representation of $\psi$ as a Legendre Transform of the Value Function of an Optimal Stopping Problem}

In this section, we first introduce an auxiliary optimal stopping problem and show that its Legendre transform is equal to the minimum probability of ruin. We also show that an individual will not purchase {\it any} deferred annuities until her wealth reaches the safe level $\bar w(A, t)$.

\subsection{A related optimal stopping problem}

Consider the following ``penalty function'' $u$ defined for $(y, t) \in [0, \infty) \times [0, T]$, for a fixed $A \in [0, c]$, by

\begin{equation}\label{3.13}
u(y, t) \triangleq \min \left( 1, \bar w(A, t) y \right).
\end{equation}

\noi Fix values $0 < \bar y < {1 \over \bar w(A, t)} < y_0$.  Note that $u$ is maximal among those functions $f$ defined on $[0, \infty) \times [0, T]$ that are concave in $y$ and satisfy

\begin{equation}\label{3.101}
f(\bar y, t) = \bar w(A, t) \bar y, \hbox{ and } f_y(\bar y, t) = \bar w(A, t),
\end{equation}

\noi and

\begin{equation}\label{3.91}
f(y_0, t) = 1, \hbox{ and } f_y(y_0, t) = 0.
\end{equation}

Define a stochastic process $Y$ by

\begin{equation}\label{3.14}
dY_s = (\ls -r) Y_s \, ds + {\mu - r \over \sigma} \, Y_s \, d
\hat B_s, \quad Y_t = y > 0,
\end{equation}

\noindent in which $\hat B$ is a standard Brownian motion with respect to a filtration $\mathbb{\hat F} = \{{\cal \hat F}_s\}_{s \geq 0 }$ of a probability space $(\hat \Omega, {\cal \hat F}, \hat \P)$, and consider the optimal stopping problem given by

\begin{equation}\label{OptStop}
\hat \psi(y, t) = \inf_{\tau \in {\cal S}_{t, T}} \hat \E \left[ \int_t^{\tau} c \, e^{-\ls (s - t)} \, Y_s \, ds + e^{-\ls (\tau - t)} \left\{ {\bf 1}_{\{\tau < T\}}  \, u(Y_\tau, \tau) + {\bf 1}_{\{\tau = T\}} \, g(Y_T) \right\} \bigg| {\cal \hat F}_t
\right],
\end{equation}

\noi in which ${\cal S}_{t, T}$ is the collection stopping times
with respect to $(\hat \Omega, {\cal \hat F}, \hat \P,
\mathbb{\hat F})$ such that $\tau$ takes values in $[t, T]$, and
in which $g$ is given by

\begin{equation}\label{terminal}
g(y) = {c-A \over r} \, y - (d - 1)  \left( {c - A \over rd} \, y
\right)^{{d \over d -1}}.
\end{equation}

\noi Recall that $d$ is given in (\ref{2.4}).

One can think of this optimal stopping problem as awarding a ``player" the running penalty $c Y_s$ between time $t$ and the time of stopping $\tau$, discounted by the probability that the player survives to time $s$.  At the time of stopping, the player receives the penalty $u(Y_\tau, \tau)$ or $g(Y_T)$ depending on whether $\tau < T$ or $\tau = T$, respectively, if she is alive then.  Thus, the player has to decide whether it is better to continue receiving the running penalty $c Y_s$ or to stop and take the final penalty $u(Y_\tau, \tau)$ or $g(Y_T)$.

Note that $\hat \psi$ is concave with respect to $y$.  Indeed, because $Y$ in (\ref{3.14}) is given by $Y_s = y \, H_s$ with
\begin{equation}
H_s = \exp \left(-(r + m - \ls)(s - t) + {\mu - r \over \sigma}
(\hat B_s - \hat B_t)  \right),
\end{equation}
\noi the integral in (\ref{OptStop}) can be written as $y \int_t^{\tau \wedge T} c \, e^{-\ls (s - t)} \, H_s \, ds$.  Thus,
because $u$ and $g$ are concave with respect to $y$, the expression in the expectation is concave with respect to $y$.  It follows that the infimum over stopping times $\tau \in {\cal S}_{t, T}$ is concave with respect to $y$.

Define the {\it continuation region} by

\begin{equation}
D = \{ (y, t) \in \R^+ \times [0, T]: \hat \psi(y, t) < u(y, t) \}
\end{equation}

\noi and consider its sections

\begin{equation}
D_t = \{ y \in \R^+: \hat \psi(y, t) < u(y, t) \}, \quad t \in [0,
T].
\end{equation}

\noi By following the same line of arguments in Section 2.7 of \cite{kn:karat2}, one can show that there are numbers $y_0(t) \ge \bar y(t) \ge 0$ (to be determined) such that $D_t = ( \bar y(t), y_0(t))$ (thus, $D = \{(y, t) \in \R^+ \times [0, T]: \bar y(t) < y < y_0(t) \}$) and that $\hat \psi$ is the unique classical solution of the following free-boundary problem (FBP):

\begin{equation}\label{FBP}
\begin{cases}
\ls \hat \psi = \hat \psi_t + (\ls - r) y \hat \psi_y + m y^2 \hat \psi_{yy} + c y \; \hbox{ on } \; D, \\
\hat \psi(y_0(t), t) = 1 \; \hbox{ and } \; \hat \psi(\bar y(t), t) = \bar w(A, t) \bar y(t) \; \hbox{ for } \; 0 \le t \le T, \\
\hat \psi(y, T) = {c-A \over r} \, y - (d -1) \left( {c - A \over rd} \, y \right)^{{d \over d -1}} \; \hbox{ for } \; y \ge 0.
\end{cases}
\end{equation}

\noi Moreover, $\hat \psi$ is $C^1$ across the free boundaries, that is,
\begin{equation}
 \hat \psi_y(y_0(t), t) = 0 \hbox{ and } \, \hat \psi_y(\bar y(t), t) = \bar w(A, t),
\end{equation}
which can be shown by using similar arguments used in the proof of Lemma 2.7.8 in \cite{kn:karat2}.

In the next section, we show that the solution of the FBP (\ref{FBP}) is intimately connected with the minimum probability of ruin.

\subsection{Relation between the FBP (\ref{FBP}) and the Minimum Probability of Ruin}

In this section, we show that the Legendre transform (see, for example, \cite{kn:karat2}) of the solution of the FBP (\ref{FBP}) is in fact the minimum probability of ruin $\psi$. To this end, note that because $\hat \psi$ in (\ref{OptStop}) is concave, we can define its convex dual via the Legendre transform. We will omit the dependence on $A$ when it is not necessary to emphasize the dependence on this variable.
\begin{thm}
Let 
\begin{equation}\label{Leg}
\Psi(w, t) \triangleq \max_{y \ge 0} [\hat \psi(y, t) - wy].
\end{equation}
for $w \ge 0$ and $t \in [0, T]$, in which 
$\hat{\psi}$ is the value function of the optimal stopping problem in (\ref{OptStop}). If 
\begin{equation}\label{eq:ass-2}
\bar{a}(s)\Psi_w(w,A,s)-\Psi_A(w,A,s)\leq 0,
\end{equation}
for  $A \in [0,c]$, $w \geq 0$ and $s \geq t$, or equivalently 
\begin{equation}\label{ineq}
\hat \psi_A(y, A, t)  \ge - y \, \bar a(t),
\end{equation}
\noi for all $y \ge 0$, $A \in [0, c]$, and $t \in [0, T]$, then $\Psi(w,t)=\psi(w,t)$.
\end{thm}

\begin{proof}
 For a given value of $t \in [0, T]$, the optimizer $y^*$ of (\ref{Leg}) solves the equation $\hat \psi_y(y, t) - w = 0$; thus, $y^*(t) = I(w, t)$, in which $I$ is the inverse of $\hat \psi_y$ with respect to $y$.  It follows that
\begin{equation}\label{3.2}
\Psi(w, t) = \hat \psi[I(w, t), t] - w I(w, t).
\end{equation}

\noi Expression (\ref{3.2}) implies that

\begin{equation}\label{3.3}
\begin{split}
\Psi_w(w, t) &= \hat \psi_y[I(w, t), t] I_w(w, t) - I(w, t) - w I_w(w, t) \\ &= w I_w(w, t) - I(w, t) - w I_w(w, t) \\ &= - I(w,
t).
\end{split}
\end{equation}
\noi Thus, $y^*(t) = I(w, t) = - \Psi_w(w, t)$. Note that from (\ref{3.3}), we have
\begin{equation}\label{3.5}
\Psi_{ww}(w, t) = -I_w(w, t) = -1/\hat \psi_{yy}[I(w, t), t],
\end{equation}
\noi and from (\ref{3.2}), we have
\begin{equation}\label{3.6}
\begin{split}
\Psi_t(w, t) &= \hat \psi_y[I(w, t), t] I_t(w, t) + \hat \psi_t[I(w, t), t] - w I_t(w, t) \\
&= w I_t(w, t) + \hat \psi_t[I(w, t), t] - w I_t(w, t) \\
&= \hat \psi_t[I(w, t), t].
\end{split}
\end{equation}

We proceed to find the boundary-value problem (BVP) $\Psi$ solves given that $\hat \psi$ solves (\ref{FBP}). In the partial differential equation (PDE) for $\hat \psi$ in (\ref{FBP}), let $y = I(w, t) = -\Psi_w(w, t)$ to obtain
\begin{equation}
\ls \hat \psi[I(w, t), t] = \hat \psi_t[I(w, t), t] + (\ls - r) I(w, t) \hat \psi_y[I(w, t), t] + m I^2(w, t) \hat \psi_{yy}[I(w,
t), t] + c I(w, t).
\end{equation}
\noi Rewrite this equation in terms of $\Psi$ to get
\begin{equation}\label{fbp1}
\ls[ \Psi(w, t) - w \Psi_w(w, t)] = \Psi_t(w, t) - (\ls - r) \Psi_w(w, t) w + m {\Psi_w^2(w, t) \over - \Psi_{ww}(w, t)} - c
\Psi_w(w, t),
\end{equation}
\noi or equivalently,
\begin{equation}\label{3.7}
\ls \Psi(w, t) = \Psi_t(w, t) + (rw - c) \Psi_w(w, t) - m{\Psi_w^2(w, t) \over \Psi_{ww}(w, t)}.
\end{equation}

Next, consider the boundary and terminal conditions in (\ref{FBP}).  First, the boundary conditions at $\bar y(t)$,
namely $\hat \psi(\bar y(t), t) = \bar w(A, t) \bar y(t)$ and $\hat \psi_y(\bar y(t), t) = \bar w(A, t)$, imply that the
corresponding dual value of $w$ is $\bar w(A, t)$ and that
\begin{equation}\label{fbp3}
\Psi(\bar w(A, t), t) = 0.
\end{equation}
  Similarly, the boundary conditions at $y = y_0(t)$, namely $\hat \psi(y_0(t), t) = 1$ and $\hat \psi_y(y_0(t), t) = 0$, imply that the corresponding dual value of $w$ is $0$ and that
\begin{equation}\label{fbp4}
\Psi(0, t) = 1.
\end{equation}
Finally, if we compute the Legendre transform of the terminal condition in (\ref{FBP}), namely $\hat \psi(y, T) =  {c-A \over r} \, y - (d -1)  \left( {c - A \over rd} \, y \right)^{{d \over d -1}}$, via the expression in (\ref{Leg}), then we obtain
\begin{equation}\label{fbp5}
\Psi(w, A, T) =  \left ( {c - A - rw \over c - A} \right)^d.
\end{equation}

Thus, we have shown that the Legendre transform $\Psi$ of the solution of the optimal stopping problem $\hat \psi$ (\ref{OptStop}), or equivalently of the FBP (\ref{FBP}), is the solution of the BVP (\ref{3.7})-(\ref{fbp5}).  
Next, we want to show that $\Psi$ equals the minimum probability of ruin $\psi$ defined in (\ref{eq:probruin}).  To this end, consider the following deferred annuitization strategy $\Ac$:  The individual buys no (additional) deferred annuity income until her wealth reaches the safe level; at that time, the individual buys a deferred life annuity to cover her remaining consumption $c - A$ and has enough wealth to sustain her until time $T$ by investing her wealth in the riskless asset.  That is, if $A_t = A$, then $\Ac$ is given by

\begin{equation}\label{eq:opt-an-pur}
A_s = A + (c-A)1_{\{s= \tau_{\bar{w}}\}},\,s\in [t,T),
\end{equation}

\noi in which

\begin{equation}
\tau_{\bar{w}} \triangleq \inf\{s \in [t,T): W_s \ge \bar{w}(A_s, s)\}.
\end{equation}

\noi Now, because $\Psi$ is the solution of the BVP (\ref{3.7})-(\ref{fbp5}), it follows that
\begin{equation}\label{eq:Psi-o-o-p}
\Psi(w, A, t) = \inf_{\pi} \P \left[ \tau_0 < \tau_d \,\big| A_t=A, W_t=w \right],
\end{equation}
for $0 \le t \le T$, in which the infimum is taken over all admissible investment strategies $\pi$ and in which the individual follows the deferred annuitization strategy $\Ac$ defined above.

Since $\Psi$ is the solution of (\ref{3.7})-(\ref{fbp5}), it follows from Lemma~\ref{lem:verf-lemma} that when assumption (\ref{eq:ass-2}) is satisfied we have that 
\begin{equation}
\Psi(w,A,t) \leq \psi(w,A,t) \quad \text{for all $w \geq 0$ and $A
\in [0,c]$.}
\end{equation}
But the last inequality together with (\ref{eq:Psi-o-o-p}) would imply that
\begin{equation}
\Psi(w,A,t) = \psi(w,A,t) \quad \text{for all $w \geq 0$ and $A
\in [0,c]$,}
\end{equation}
and that an individual will not purchase {\it any} deferred annuities until her wealth reaches the safe level $\bar w(A, t)$ (see (\ref{eq:opt-an-pur})). It is clear that (\ref{ineq}) is equivalent to (\ref{eq:ass-2}).
\end{proof}

To compute the solution of this BVP, one can apply the Projected SOR method (see, for example, \cite{howison})
to compute the solution of the (\ref{FBP}), then apply the expression in (\ref{Leg}) to get $\Psi$.

\section{Numerical Example}

\begin{figure}[h]
\begin{center}

\includegraphics[width = 0.9\textwidth,height=6.5cm]{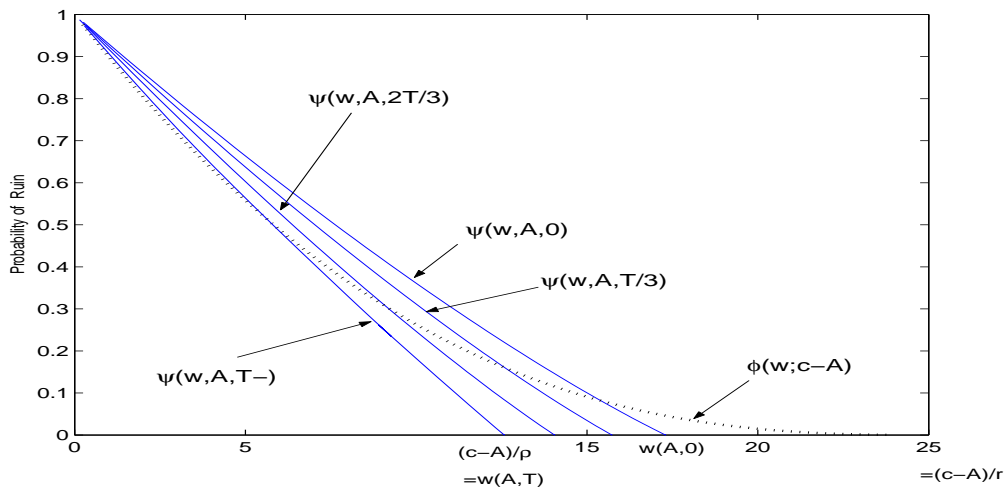}
\caption{The minimum probability of ruin in the presence of
deferred annuities whose benefits will be active at time $T$.}
\label{fig:I}
\end{center}
\end{figure}

In this section, we present a numerical example to demonstrate our
results. The FBP as given is (\ref{FBP}) is not readily amenable
to numerical solution because the free-boundary is unknown.
However, solving the optimal stopping problem (\ref{OptStop}), or
equivalently solving the FBP (\ref{FBP}), it is sufficient to
solve the following variational inequality:
\begin{equation}\label{VI}
\begin{cases}
\max \left[ \ls \hat \psi - \hat \psi_t - (\ls - r)y \hat \psi_y -
m y^2 \hat \psi_{yy} - c y, \hat \psi - u \right] = 0, \\ \hat
\psi(y, T) = {c-A \over r} \, y - (d - 1)  \left( {c - A \over rd}
\, y \right)^{{d \over d -1}}.
\end{cases}
\end{equation}
We numerically solve this variational inequality using the
projected SOR method. For our numerical example we choose the
following values for the parameters of the problem:
\begin{itemize}
\item $\lambda^{S}=\lambda^{O}=0.02$; the expected future lifetime
is 50 years,
\item the riskless interest rate $r$ is $2\%$ over inflation,
\item the appreciation rate of the risky investment is $\mu=6\%$
over inflation,
\item the volatility of the risky asset is $20 \%$,
\item the individual consumes $c=1.5$ units of wealth per year,
\item annuity income from previous investments is $A=1$ units of
wealth per year,
\item the deferred annuity benefits start in $T=5$ years.
\end{itemize}
In Figure~\ref{fig:I} we graph the probability of ruin with and without deferred annuities, $\psi(w,A,t)$ and $\phi(w;c-A)$,
respectively.  We confirm the convexity of $\psi(w,A,t)$ with respect to wealth $w \in [0, \bar w(A, t)]$.  Next, we observe that the presence of annuities reduces the safe level from $(c-A)/r$ to $\bar{w}(A,t)$.  We see that the safe level $\bar w(A,t)$ decreases with respect to time $t \in [0, T]$, which can be confirmed by taking the derivative of (\ref{eq:bar-w}).  Next, we observe that the probability of ruin $\psi(w,A,t)$ decreases with respect to time $t \in [0, T]$. We also observe that the optimal investment in the risky asset
\begin{equation}
w \rightarrow
\pi^*(w,A,t)=-\frac{\mu-r}{\sigma^2}\frac{\Psi_w(w,A, t)}{\Psi_{ww}(w,A,t)},
\end{equation}
is not continuous at $w=\bar{w}(A,t)$, since $\psi_w(w,A, t)<0$ for $w<\bar{w}(A,t)$ and $\psi(w,A, t)=0$ for $w>w(A,t)$.

It is interesting to note that $\phi(w; c-A)$ is less than $\psi(w, A, t)$, for wealth low enough when $t < T$, because $\phi(w;c-A)$ is calculated assuming that the investor receives the benefit $A$ {\it immediately}, whereas $\psi(w,A,t)$ is calculated assuming that she will receive the benefit $A$ at time $T$.  In other words, an individual with low wealth prefers immediate income with {\it no} opportunity to buy deferred annuities, instead of deferred annuity income together with an expanded opportunity set for buying further deferred income.

We also numerically validated (\ref{ineq}), from which we conclude that the investor buys a deferred life annuity only if she has enough money to buy one to cover her excess consumption and to consume the remainder which is invested in the riskless asset until time $T$ without ruining (see Remark 3.1 and Section 4.2). Then at time $T$, the deferred annuity will cover all her consumption and she will avoid ruin. The qualitative properties of the probability of ruin problem that we see in this example appear to be true for every possible range of parameters.

\section{Summary and Conclusions}

We showed that an individual will buy a deferred annuity only when
her wealth reaches the safe level $\bar w(A, t)$ given in
(\ref{eq:bar-w}).  When wealth is lower than that amount, the
individual does not want to relinquish the ability to trade on the
volatility of the risky asset, as we observed in the simple
example in Section 2.2.

One could add immediate life annuities into the financial market,
as in \cite{mmy}, with the understanding that the minimum
probability of ruin in this case, say $\tilde \psi = \tilde
\psi(w, A, t)$, is a function of wealth $w$ at time $t$ and
deferred annuity income $A$ purchased on or before time $t$ to
begin at time $T$.  Recall that we allow $A$ to include pension or
Social Security income that will begin at time $T$.  In this case,
one would obtain a result similar to what we obtained in the case
of deferred life annuities.  Specifically, the individual would
not purchase immediate annuities before time $T$ unless wealth
were to reach the safe level $\bar {\bar w}(A, t)$ given by

\begin{equation}\label{safe:immed}
\bar {\bar w}(A, t) = \min \left[ A \frac{1-e^{-r(T - t)}}{r} +
(c-A) {1 \over \rho}, \; {c \over \rho} \right], \quad t \leq T,
\, A \in [0,c],
\end{equation}

\noi which is less than $\bar w(A, t)$ because $\rho > r$.
Therefore, in the presence of immediate annuities, an individual
would {\it never} purchase deferred annuities.  Perhaps our
results confirm the low level of participation in  deferred
annuity markets.

\bibliographystyle{dcu}
\bibliography{references}

\end{document}